\documentclass[12pt,reqno]{amsart}

\usepackage{amscd,amssymb,amsmath}

\usepackage{tikz-cd}

\makeatletter
\def\@tocline#1#2#3#4#5#6#7{\relax
 \ifnum #1>\c@tocdepth 
  \else
    \par \addpenalty\@secpenalty\addvspace{#2}%
   \begingroup \hyphenpenalty\@M
    \@ifempty{#4}{%
      \@tempdima\csname r@tocindent\number#1\endcsname\relax
    }{%
      \@tempdima#4\relax
    }%
    \parindent\z@ \leftskip#3\relax \advance\leftskip\@tempdima\relax
    \rightskip\@pnumwidth plus4em \parfillskip-\@pnumwidth
    #5\leavevmode\hskip-\@tempdima
      \ifcase #1
       \or\or \hskip 1em \or \hskip 2em \else \hskip 3em \fi%
      #6\nobreak\relax
    \dotfill\hbox to\@pnumwidth{\@tocpagenum{#7}}\par
    \nobreak
    \endgroup
  \fi}
 \makeatother
 
\newcommand{\bZ}{\mathbb Z}

\newcommand{\bC}{\mathbb C}
\newcommand{\bR}{\mathbb R}

\DeclareMathOperator{\Id}{Id}

\DeclareMathOperator{\Span}{span}
\DeclareMathOperator{\Lef}{Lef}
\DeclareMathOperator{\Tr}{Tr}
\DeclareMathOperator{\SW}{SW}
\DeclareMathOperator{\FO}{FO}
\DeclareMathOperator{\red}{red}

\theoremstyle{plain}

\newtheorem*{thmA}{Theorem A}

\theoremstyle{definition}

\newtheorem{ex}{Example}

\theoremstyle{remark}

\begin{document}

\title[Splitting Formula]{A Splitting Formula in Instanton Floer Homology}
\author{Nima Anvari}
\date{\today}
\email{anvarin@math.mcmaster.ca}
\thanks{The author thanks Nikolai Saveliev for many helpful discussions, valuable suggestions and  for encouragement in writing this paper.}

\begin{abstract}
In a recent paper, Lin, Ruberman and Saveliev proved a splitting formula expressing the Seiberg-Witten invariant $\lambda_{\SW}(X)$ of a smooth $4$-manifold with rational homology of $S^1\times S^3$ in terms of the Fr{\o}yshov invariant $h(X)$ and a Lefschetz number in reduced monopole Floer homology. In this note we observe that a similar splitting formula holds in reduced instanton Floer homology. 
\end{abstract}

\maketitle


\section{Introduction}

Let $X$ be a smooth, oriented spin 4-manifold that has integral homology of $S^1 \times S^3$. A smooth invariant of $X$, denoted by $\lambda_{\SW}(X)$, was defined by Mrowka, Ruberman and Saveliev  \cite{Mrowka_2011} as a signed count of Seiberg--Witten monopoles plus an index-theoretic correction term. Recently, Lin, Ruberman and Saveliev \cite{lin2018splitting} extended the definition of $\lambda_{\SW} (X)$ to all smooth, oriented spin 4-manifolds with rational homology of $S^1 \times S^3$ and proved a splitting formula relating $\lambda_{\SW}(X)$ to the Fr{\o}yshov invariant $h(Y,\mathfrak{s})$, where $Y$ is an embedded rational homology $3$-sphere generating $H_3(X;\bZ)$ (assuming such exists) and $\mathfrak{s}$ is the induced spin structure. It is proved in \cite{froyshov2010monopole} that $h(Y,\mathfrak{s})$ is an invariant of $X$, denoted by $h(X)$. The splitting formula relating these two invariants reads \cite[Theorem A]{lin2018splitting}
\[
\lambda_{SW}(X) + h(X) = - \Lef ( W_{\ast}: HM^{\red}(Y,\mathfrak{s} )\rightarrow HM^{\red}(Y,\mathfrak{s})),
\]
where $W$ is the spin cobordism from $Y$ to itself obtained by cutting $X$ open along $Y$ and $HM^{\red}(Y,\mathfrak{s})$ is the reduced monopole Floer homology. 

In this note we show that a similar splitting formula holds  when reduced monopole Floer homology is replaced by reduced instanton Floer homology, $\lambda_{\SW}(X)$ is replaced by the Furuta--Ohta invariant $\lambda_{\FO}(X)$, as defined in \cite{FURUTA1993291}, and the monopole Fr{\o}yshov invariant $h(Y)$ is replaced by the instanton Fr{\o}shov invariant $h(Y)$ (that the same notation is used for both Fr{\o}yshov invariants may be confusing, although the two invariants are known to coincide in all calculated examples).

More precisely, let $X$ be a $\bZ[\bZ]$-homology $S^1 \times S^3$, that is, a smooth 4-manifold such $H_{\ast}(X; \bZ) = H_{\ast}(S^1 \times S^3; \bZ)$ and $H_{\ast}(\widetilde{X};\bZ) = H_{\ast}(S^3;\bZ)$, where $\widetilde{X}$ is the universal abelian cover of $X$. The moduli space $\mathcal{M}^{\ast}(X)$ of irreducible anti-self dual connections  on a trivial $SU(2)$ bundle over $X$ is compact and has formal dimension zero; in fact, it coincides with the moduli space of flat $SU(2)$ connections on $X$. With an appropriate choice of orientations and admissible perturbations, the Furuta--Ohta invariant  is defined as a signed count 
\[
\lambda_{\FO}(X) = \frac 1 4\, \#  \mathcal{M}^{\ast}(X),
\]
see Ruberman--Saveliev \cite{Ruberman_2004} for details.

\begin{thmA}\label{T:main}
Let $X$ be a $\bZ[\bZ]$-homology $S^1\times S^3$ and suppose that there exists an embedded integral homology sphere $Y$ which generates $H_3(X;\bZ)$. Cut $X$ open along $Y$ to obtain a homology cobordism $W$ from $Y$ to itself, and denote by $W_{\ast}$ and $\widehat{W}_{\ast}$ the induced homomorphisms in, respectively, the instanton Floer homology and the reduced instanton Floer homology of $Y$. Then the quantity
\begin{equation}\label{E:h(X)}
h(X) =  \dfrac{1}{2}\, ( \Lef(\widehat{W}_{\ast})- \Lef(W_{\ast}))
\end{equation}
defined in terms of the Lefschetz numbers is independent of $W$ and coincides with the instanton Fr{\o}yshov invariant $h(Y)$. Moreover, the following splitting formula holds
\begin{equation}\label{E:splitting}
\lambda_{\FO}(X) +  h(X) = \dfrac{1}{2}\,\Lef (\widehat{W}_{\ast}: \widehat{HF^{\ast}}(Y) \longrightarrow  \widehat{HF^{\ast}}(Y)).
\end{equation}
\end{thmA}

We refer the reader to Fr{\o}yshov \cite{FROYSHOV2002525} for the definitions of the reduced instanton Floer homology and the instanton Fr{\o}yshov invariant; see also Section \ref{S:floer} and Section \ref{S:proof}. Note that the splitting formula \eqref{E:splitting} follows easily from \eqref{E:h(X)} and the formula $\lambda_{\FO}(X) = 1/2\,\Lef(W_{\ast})$ of \cite{Ruberman_2004}. To prove that $h(X)$ is a well-defined invariant, we will show that it coincides with the instanton Fr{\o}yshov invariant $h(Y)$ for any choice of $Y$. When $W$ is a product cobordism $ Y\times [0,1]$, both homomorphisms $W_{\ast}$ and $\widehat{W}_{\ast}$ are identity maps, and \eqref{E:h(X)} reduces to the Fr{\o}yshov formula 
\[
h(X) = h(Y) = \frac{1}{2}(\chi(\widehat{HF_{\ast}}(Y)) - \chi(HF_{\ast}(Y))).
\]
We will prove the general case of \eqref{E:h(X)} in this paper using properties of the special boundary maps of Fr{\o}yshov \cite{FROYSHOV2002525}.

\begin{ex}\label{Ex:one}
Let $Y$ be the Brieskorn homology $3$-sphere $\Sigma(2,7,13)$ oriented as the link of a complex surface singularity, and $\tau$ the involution on $\Sigma(2,7,13)$ induced by complex conjugation. The mapping torus of $\tau$ is a smooth $4$-manifold $X_{\tau} = [0,1] \times Y/ (0,x) \sim (1,\tau(x))$ and a $\bZ[\bZ]$-homology $S^1 \times S^3$. The mod $8$ graded instanton Floer homology of $Y$ is given by (see \cite{fintushel1990instanton})
\[
HF_{\ast}(Y) = (0, \bZ^4,0,\bZ^2, 0, \bZ^4,0,\bZ^2)
\]
and the reduced instanton Floer homology 
\[
\widehat{HF}_{\ast}(Y) = (0, \bZ^2,0,\bZ^2, 0, \bZ^2,0,\bZ^2)
\]
is completely determined by $HF_{\ast}(Y)$ and the invariant $h(Y)$ which was computed by Fr{\o}yshov \cite{froyshov2004inequality} to be $h(Y) = 2$. In this case the Furuta--Ohta invariant $\lambda_{\FO} (X_{\tau})$ equals the Neumann--Siebenmann $\overline{\mu}$--invariant; see \cite{Ruberman2004}, and a closed form expression is given by Saveliev \cite[Theorem 6.28, p.147]{Nikolai} as $-b_1+b_3$, where $b_i$ is the rank of $HF_i(Y)$. Therefore, $\lambda_{\FO}(X_{\tau}) = -2$. On the other hand, let $W$ be the homology cobordism obtained by cutting $X_{\tau}$ open along $Y$. According to \cite[Proposition 9.2]{Ruberman2004}, the induced map $W_{\ast} = \tau_{\ast}: HF_{k}(Y) \rightarrow HF_k(Y)$ is the identity for $k \equiv 1 \pmod 4$ and minus the identity for $k \equiv -1 \pmod 4$.  It follows that $\Lef(W_{\ast})=-4$ and $\Lef(\widehat{W}_{\ast})=0$ hence
\[
h(X_{\tau}) = \dfrac{1}{2}\,(\Lef(\widehat{W}_{\ast}) - \Lef(W_{\ast}) ) =2, 
\]
which matches $h(Y)$ and hence confirms the splitting formula $\lambda_{\FO}(X_{\tau}) + h(X_{\tau}) = 0=1/2\, \Lef(\widehat{W}_{\ast})$.
\end{ex}

\begin{ex}\label{Ex:two}
Let $W$ be the Akbulut cork as in \cite[Section 9.3]{Ruberman_2004}. That is, $W$ is a smooth contractible manifold with boundary an integral homology sphere $\Sigma$ that can be embedded into a blown up elliptic surface $E(n)\, \#\, \overline{\bC\rm P}^2$ in such a way that cutting it out and re-gluing by an involution $\tau: \Sigma \rightarrow \Sigma$ changes the smooth structure on $E(n)\, \#\, \overline{\bC\rm P}^2$ but preserves its homeomorphism type. The instanton homology groups are
\[
HF_{\ast}(\Sigma) = (0,\bZ,0,\bZ,0,\bZ,0,\bZ)
\]
and, since $\Sigma$ is homology cobordant to zero, $h(\Sigma)=0$ and $\widehat{HF}_{\ast}(\Sigma)=HF_{\ast}(\Sigma)$. Let $X_{\tau}$ be the mapping torus of $\tau$ then it follows from \cite{Ruberman_2004} that $W_{\ast} = -\Id$ and $\lambda_{\FO}(X_{\tau}) = 2$, which again confirms the splitting formula.
\end{ex}


\section{Reduced instanton Floer (co)homology}\label{S:floer}
In this section, we recall the definition of the reduced instanton Floer homology (cohomology) groups; the reader is referred to Fr{\o}yshov \cite{FROYSHOV2002525} for all the details. We work throughout with real coefficients and the orientation conventions of \cite{FROYSHOV2002525} and use the canonical identification $HF_q (Y) = HF^{5-q} (\overline{Y})$. 

The instanton Floer cohomology $HF^*(Y)$ is defined as homology of the mod 8 graded instanton cochain complex $CF^*(Y)$ generated by the irreducible flat $SU(2)$ connections on $Y$, with respect to the boundary map $d: CF^*(Y) \to CF^{*+1} (Y)$ given by a signed count of anti-self dual connections on $\mathbb R \times Y$ of finite energy. The definition of the reduced instanton cohomology further employs special boundary maps $\delta_0$ and $\delta_0^{\prime}$, which are defined as follows. Let $\theta$ be the trivial flat connection on a trivialized $SU(2)$ bundle over $Y$. Define 
\[
\delta:\; CF^4(Y) \longrightarrow  \bR\quad  \text{by}\quad
\delta(\alpha) = \#\mathcal{\check {M}}(\theta,\alpha),
\]
where $\mathcal{M}(\theta,\alpha)$ is the moduli space of anti-self dual connections on $\bR \times Y$ limiting to $\theta$ at $+\infty$ and $\alpha$ at $-\infty$, and $\check {\mathcal M} (\theta,\alpha)$ its quotient by translations. Define 
\[
\delta^{\prime}: \bR \longrightarrow CF^1(Y)\quad\text{by}\quad \delta^{\prime}(1) = \sum_{\beta}\#\check{\mathcal M}(\theta, \beta)\cdot \beta,
\]
with the summation extending over the generators $\beta \in CF^1(Y)$. The maps $\delta$ and $\delta^{\prime}$ satisfy equations $d\delta=0$ and $d\delta^{\prime}=0$, thereby inducing maps in cohomology,
\[
\delta_{0}: HF^4(Y) \longrightarrow \bR\quad \text{and}\quad \delta^{\prime}_0: \bR \longrightarrow HF^1(Y).
\]
The special boundary maps $\delta_0$ and $\delta_0^{\prime}$ are further included into the sequence of maps
\[
\delta_n = [\delta v^n]: HF^{4-4n}(Y) \longrightarrow \bR
\quad\text{and}\quad 
\delta_n^{\prime} = [v^n \delta^{\prime}] : \bR \longrightarrow HF^{1+4n}(Y)
\]
using the homomorphism $v: CF^{\ast} (Y) \longrightarrow CF^{\ast+4}(Y)$ defined in \cite{FROYSHOV2002525}. It follows from the cochain homotopy formula \cite[Theorem 4(ii)]{FROYSHOV2002525} that either $\delta_0$ or $\delta_0^{\prime}$ must vanish; moreover, if $\delta_0 = 0$ then $\delta_n=0$ for all $n$, and similarly for $\delta_n^{\prime}$.

We also need to recall the relation between the special boundary maps and the homomorphisms in Floer homology induced by negative definite cobordisms. Let $W$ be a connected Riemannian $4$-manifold with two cylindrical ends, $\bR_{-} \times Y_1$ and $\bR_{+} \times Y_2$, and assume that $W$ has negative definite intersection form, $H_1 (W;\bZ) = 0$, and both $Y_1$ and $Y_2$ are integral homology spheres. Then $W$ induces a degree preserving cochain homomorphism 
\medskip
\[
W^{\ast}: CF^{\ast}(Y_1) \rightarrow CF^{\ast}(Y_2)\quad \text{by} \quad W^*(\alpha)\;=\; 
\sum_{\beta} \#\mathcal{M}(W; \beta,\alpha)\cdot \beta,
\]
where the summation extends over the generators $\beta \in CF^{\ast}(Y_2)$ with the same index as $\alpha$, and $\mathcal M(W; \beta,\alpha)$ is the zero-dimensional part of the moduli space of anti-self dual connections on $W$ limiting to $\alpha$ at $-\infty$ and to $\beta$ at $+\infty$. Fr{\o}yshov \cite[Theorem 7]{FROYSHOV2002525} shows that $\delta_0 W^{\ast} = \delta_0$ and $W^{\ast}\delta_0^{\prime}= \delta_0^{\prime}$ and, moreover, that there exist integers $a_{ij}$ and $b_{ij}$ such that 
\begin{align}
\delta_n W^{\ast} = \delta_n + \sum_{i=0}^{n-1}a_{in}\delta_{i}, \label{E:delta}\\
W^{\ast}\delta^{\prime}_n = \delta_n^{\prime} + \sum_{i=0}^{n-1}b_{in}\delta_i^{\prime}.\label{E:delta'}
\end{align}
A straightforward grading count shows that $a_{in} = 0$ and $b_{in} = 0$ whenever $i$ and $n$ have opposite parity.

We are now ready to define the reduced instanton cohomology groups. Let $B^{\ast} \subseteq  HF^{\ast} (Y)$ denote the linear span of the vectors $\delta_n^{\prime} (1)$ for all $n$ so that $B^1 = \Span\{{\delta_{2k}^{\prime}}(1)\}$ and $B^5 = \Span\{{\delta_{2k+1}^{\prime}}(1)\}$ with $k \ge 0$, while $B^q = 0$ for $q \ne 1, 5 \pmod 8$. In addition, let 
\[
Z^{\ast} = \bigcap_n\; \ker(\delta_n)\; \subseteq\; HF^{\ast} (Y)
\]
so that 
\[
Z^0 = \bigcap_{k \ge 0}\;\ker(\delta_{2k+1}) \subseteq HF^0(Y)\quad\text{and}\quad Z^4 = \bigcap _{k \ge 0}\; \ker(\delta_{2k}) \subseteq HF^4(Y),
\]
while $Z^q = HF^q(Y)$ for $q \neq 0, 4\pmod 8$. The reduced Floer cohomology groups are then defined as
\[
\widehat{HF}^{\,q} (Y) = Z^q / B^q.
\]
For any negative definite cobordism $W$ as above, the map $W^{\ast}$ leaves the subspaces $Z^q$ and $B^q$ invariant and induces a map \cite[Theorem 8]{FROYSHOV2002525}
\[
W^{\ast}: \widehat{HF}^{\ast}(Y_1) \longrightarrow \widehat{HF}^{\ast}(Y_2).
\]

Finally, we mention that both instanton Floer cohomology $HF^{\ast} (Y)$ and reduced instanton Floer homology $\widehat{HF}^{\ast} (Y)$, which \emph{a priori} have a mod 8 grading, are 4-periodic; see \cite[Theorem 2 and Corollary 3]{FROYSHOV2002525}.


\section{Proof of Theorem A}\label{S:proof}
Recall that the (instanton) Fr{\o}yshov invariant $h(Y)$ is defined by the formula 
\[
h(Y)=\dfrac{1}{2}\,(\chi(HF^{\ast}(Y))-\chi(\widehat{HF}^{\ast}(Y)).
\]
The theorem will be proved as soon as we show that $h(X) = h(Y)$ for any choice of $Y$. Using the aforementioned 4-periodicity in Floer homology, we only have two cases to consider, depending on which one of the special boundary maps, $\delta_0$ or $\delta_0^{\prime}$, vanishes.

\medskip

\noindent\textbf{Case 1.} Suppose that $\delta_0^{\prime} = 0$ and hence $\delta_n^{\prime} =0$ for all $n$. This implies that $B^1=0$ and $\widehat{HF}^1 (Y) = Z^1 = HF^1 (Y)$. Thus the reduced theory differs from the unreduced one only in degree 0 mod 4, and 
\begin{align*}
h(X) &= \dfrac{1}{2}\,( \Lef(W^{\ast}) - \Lef(\widehat{W}^{\ast}) )  \\
&= \Tr\,(W^{\ast}\colon HF^0(Y) \rightarrow HF^0(Y)) - \Tr\,(\widehat{W}^{\ast}\colon \widehat{HF}^0(Y) \rightarrow \widehat{HF}^0(Y)) .
\end{align*}
Note that, since $B^0=0$, we have 
\[
\widehat{HF}^0(Y) = Z^0 = \bigcap_{k \ge 0}\;\ker(\delta_{2k+1}).
\]
We claim that $h(X) = \dim(HF^0(Y)/Z^0)$; this will imply that $h(X)$ equals the Fr{\o}yshov invariant $h(Y)=\dim HF^0(Y) - \dim \widehat{HF}^0(Y)$.


To prove the claim, we first observe that, since $\delta_1 W^{\ast} = \delta_1$, we have a commutative diagram of exact sequences
\medskip
\begin{center}
\begin{tikzcd}
    0\arrow{r} & Z^0_1 \arrow{r} \arrow{d}{\widehat{W}_1^{\ast}} &  HF^0(Y)\arrow{r}{\delta_1} \arrow{d}{W^{\ast}} & \bR \arrow{d}{\Id}  \\
    0\arrow{r} & Z^0_1 \arrow{r} &  HF^0(Y) \arrow{r}{\delta_1} & \bR
\end{tikzcd}
\end{center}

\medskip\noindent
where $Z^0_1 = \ker(\delta_1)$ and $\widehat{W}^{\ast}_1$ is the restriction of $W^{\ast}$ to $Z^0_1$. It follows that $\Tr(W^{\ast}) - \Tr(\widehat{W}_1^{\ast}) = 1$ if $\delta_1 \neq 0$ and $\Tr(W^{\ast}) - \Tr(\widehat{W}_1^{\ast}) = 0$ if $\delta_1 = 0$. Next we use relation \eqref{E:delta} to define a sequence of $W^{\ast}$--invariant subspaces
\[
Z^0_k\; =\; \ker(\delta_1)\, \cap \ldots  \cap\, \ker(\delta_k)\; \subseteq\; HF^0(Y)
\]
for $k$ odd, and construct the following commutative diagrams with exact rows,
\begin{center}
\begin{tikzcd}
    0\arrow{r} & Z^0_k  \arrow{r} \arrow{d}{\widehat{W}_k^{\ast}} & Z^0_{k-2 }\arrow{r}{\delta_k} \arrow{d}{\widehat{W}_{k-2}^{\ast}} & \bR \arrow{d}{\Id} \\
    0\arrow{r} & Z^0_k \arrow{r} &  Z^0_{k-2}  \arrow{r}{\delta_k} & \bR
\end{tikzcd}
\end{center}

\medskip\noindent
where $\widehat{W}^*_k$ denotes the restriction of $W^*$ to subspace $Z^0_k$. It follows that $\Tr(\widehat{W}^{\ast}_{k-2}) - \Tr(\widehat{W}^{\ast}_{k}) = 1$ if $\delta_k \neq 0$ and $\Tr(\widehat{W}^{\ast}_{k-2}) - \Tr(\widehat{W}^{\ast}_{k}) = 0$ if $\delta_k = 0$. Apply induction to the sequence of subspaces
\[
Z^0_m\; \subseteq \ldots \subseteq\; Z^0_3\; \subseteq\; Z^0_1\; \subseteq\; HF^0(Y)
\]
to conclude that
\[
\Tr\,(W^{\ast})\, =\, \dim(HF^0(Y)/Z^0_m)\, +\, \Tr\,(\widehat{W}^{\ast}_m)
\]
for all odd $m$. Since $Z^0 = Z^0_m$ and $\widehat{W}^* = \widehat{W}^{\ast}_m$ for all sufficiently large odd $m$, the claim follows. 

\medskip
\noindent\textbf{Case 2}. Suppose that $\delta_0=0$ and hence $\delta_n=0$ for all $n$. In this case, $\widehat{HF^0}(Y)=HF^0(Y)$ and $\widehat{HF^1}(Y) = HF^1(Y)/B^1$. The reduced Floer cohomology differs from Floer cohomology only in degree 1 mod 4, and
\begin{align*}
h(X) &= \dfrac{1}{2}\,( \Lef(W^{\ast}) - \Lef(\widehat{W^{\ast}}) )  \\
&= \Tr\,(\widehat{W}^{\ast}\colon \widehat{HF}^1(Y) \rightarrow \widehat{HF}^1(Y)) - \Tr\,(W^{\ast}\colon HF^1(Y) \rightarrow HF^1(Y))
\end{align*}
It follows from the commutative diagram of exact sequences
\begin{center}
\begin{tikzcd}
0\arrow{r} & B^1  \arrow{r}{} \arrow{d}{W^{\ast}_1} & HF^1(Y) \arrow{r}{} \arrow{d}{W^{\ast}} & \widehat{HF^1}(Y) \arrow{r}\arrow{d}{\widehat{W}^{\ast}} & 0 \\
0\arrow{r} & B^1 \arrow{r}{} & HF^1(Y)   \arrow{r}{} &  \widehat{HF}^1(Y)  \arrow{r} & 0
\end{tikzcd}
\end{center}
that $\Tr(W^{\ast}) - \Tr(\widehat{W}^{\ast})= \Tr\,(W^{\ast}_1)$, where $W^{\ast}_1$ is the restriction of $W^{\ast}$ to $B^1$. We claim that 
\[
\Tr\,(W^{\ast}_1) = \dim(B^1)= \dim HF^1(Y) - \dim \widehat{HF}^1(Y). 
\]
This will imply that $h(X)$ equals $h(Y) = \dim \widehat{HF}^1 (Y) - \dim HF^1(Y)$ so the conclusion will follow.

To prove the claim, use relation \eqref{E:delta'} to define $W^{\ast}$--invariant subspaces $B^1_k = \Span\{ \delta^{\prime}_0(1), \dots, \delta^{\prime}_k(1) \}$ for $k$ even. We have an increasing sequence of quotient spaces,
\[
\widehat{HF}^1(Y) = HF^1(Y)/B^1_m\; \subseteq\; HF^1(Y)/B^1_{m-2}\; \subseteq \ldots 
\subseteq\; HF^1(Y)/B^1_0
\]
corresponding to the sequence $B^1_0\; \subseteq \ldots \subseteq B^1_m\;  \subset HF^1(Y)$, where $B^1 = B^1_m$ for some sufficiently large even $m$. We also have the commutative diagrams with exact rows
\medskip
\begin{center}
\begin{tikzcd}
\bR  \arrow{r}{\delta^{\prime}_k} \arrow{d}{\Id} & HF^1(Y)/B^1_{k-2}  \arrow{r}{} \arrow{d} {\overline{W}^{\ast}_{k-2}} & HF^1(Y)/B^1_{k} \arrow{r}\arrow{d}{\overline{W}^{\ast}_{k}} & 0 \\
\bR  \arrow{r}{\delta^{\prime}_k} & HF^1(Y)/B^1_{k-2}   \arrow{r}{} &  HF^1(Y)/B^1_{k}  \arrow{r} & 0
\end{tikzcd}
\end{center}

\medskip\noindent
and the relations $\Tr\,(\overline{W}^{\ast}_{k-2}) - \Tr\,(\overline{W}^{\ast}_{k}) = 1$ if $\delta^{\prime}_k \neq 0$ and $\Tr\,(\overline{W}^{\ast}_{k-2}) - \Tr\,(\overline{W}^{\ast}_{k}) = 0$ if $\delta^{\prime}_k = 0$. By induction we obtain 
\[
\Tr( \overline{W}^{\ast}_{0}) = \dim(B^1/B^1_0) + \Tr(\widehat{W}^{\ast}).
\]
Finally, it follows from the relation $W^{\ast}\delta^{\prime}_0 = \delta^{\prime}_0$ and the exact sequence
\[
0 \longrightarrow B^1_0 \longrightarrow HF^1(Y) \longrightarrow HF^1(Y)/B^1_0 \longrightarrow 0
\]
that $\Tr(\overline{W}^{\ast}_0) = \Tr(W^{\ast}) - \dim (B^1_0)$. Therefore, $\Tr(W^{\ast}) - \Tr(\widehat{W}^{\ast})= \dim(B^1)$, which completes the proof of the claim.

\bibliographystyle{alpha}
\bibliography{Bibliography}

\end{document}